\documentclass{gtart}
\usepackage{amssymb,epsf}
\input gtmonout
\volumenumber{1}
\volumeyear{1998}
\volumename{The Epstein birthday schrift}
\pagenumbers{127}{138}
\received{26 August 1997}
\published{22 October 1998}
\papernumber{6}

\newtheorem{theorem}{Theorem}[section]
\newtheorem{lemma}[theorem]{Lemma}

\newtheorem{corollary}[theorem]{Corollary}

\theoremstyle{definition}
\newtheorem{question}[theorem]{Question}
\newtheorem{remark}[theorem]{Remark}

\def\xb{\overline{x}}
\def\yb{\overline{y}}

\begin{document}
\title {On the Burau representation modulo a small prime}
\authors{D Cooper\\D\thinspace D Long}
\shortauthors{D Cooper and D\thinspace D Long}
\asciiauthors{D Cooper and D D Long}

\begin{abstract}
We discuss techniques for analysing the structure of the group obtained
by reducing the image of the Burau representation of the braid group modulo a prime.
The main tools are a certain sesquilinear form first introduced by Squier
and consideration of the action of the group on a Euclidean building.
\end{abstract}


\maketitle

\section{Introduction}
Despite the work of many authors, the group theoretic image of linear representations 
of the braid groups remains mysterious in most cases. The first nontrivial example,
the {\em Burau representation} is not at all well understood. This representation
\[ \beta_n \co B_n \rightarrow GL( n-1 , {\bf Z}[t,t^{-1}]) \]
is known not to be faithful for $n \geq 6$ (\cite{Mo} and \cite{LP}) but the nature
of the image group and in particular, a presentation for the image group has not been found.
In \cite{CL}, we simplified the problem by composing $\beta_n$ with the map which reduces
coefficients modulo $2$. In this way, we were able to give a presentation for the image
of the simplified representation $\beta_4\otimes{\bf Z}_2$. (Throughout this paper we
use the notation ${\bf Z}_p$ for the finite field with $p$ elements.) Of course, the 
motivation for this approach comes from the classical problem of whether the representation
 $\beta_4$ is faithful and to this end we pose the question:
\begin{question}
Is there any prime $p$ for which the representation
\[ \beta_4\otimes{\bf Z}_p \co B_4 \rightarrow GL( 3 , {\bf Z}_p[t,t^{-1}]) \]
is faithful?
\end{question}
It is a consequence of some results of this note that the representation
is not faithful in the case $p = 3$, (below we exhibit a braid word in the kernel) however
the program for attacking the problem runs into difficulty at the final stage when $p = 5$.
This case remains open and has some features which suggest it may be different to the
first two primes.

In order to describe our approach, we recall that the group 
$GL( 3 , {\bf Z}_p[t,t^{-1}])$ acts on a certain contractible two dimensional simplicial
complex, $\Delta = \Delta(p)$ a so-called {\em Euclidean building} (see \cite{Br}). 
This is defined by 
embedding \[ GL( 3 , {\bf Z}_p[t,t^{-1}]) \longrightarrow GL( 3 , {\bf Z}_p(t)) \]
where ${\bf Z}_p(t)$ is the field of fractions of the ring ${\bf Z}_p[t,t^{-1}]$. This
target group admits a discrete rank one valuation defined by $\nu(p/q) = degree(q) - degree(p)$.
A standard construction  now yields the complex $\Delta$. We briefly outline
how this building and action are defined, restricting our
attention to the case $n = 4$, since this is the only case in which we shall subsequently be
interested. This will serve the additional purpose of establishing notation.  Standard
properties of $\nu$ imply that  
$${\cal O} = \{\, x \in {\bf Z}_p(t)\;\; | \;\; \nu(x) \geq 0 \,\} $$
is a subring of ${\bf Z}_p(t)$, the {\em valuation ring} associated to $\nu$.
This is a local ring and the unique maximal ideal is  easily seen to be 
${\cal M} = \{\, x \in {\bf Z}_p(t)\;\; | \;\;
\nu(x) > 0 \,\} $, a principal ideal. Choose some generator $\pi$ for this  ideal. This 
element is called a {\em uniformizing parameter} and by construction we have that
$\nu(\pi) = 1 $. Since $\cal M$ is  maximal, the quotient $k = {\bf Z}_p(t)/{\cal M}$ 
is a field, the {\em residue class field}. One sees easily that in this case, the  
residue class field is ${\bf Z}_p$.

Now let $V$ be the vector space ${\bf Z}_p(t)^3$.  By a {\em lattice} in $V$ we shall 
mean an $\cal O$--submodule, $L$, of the form $L = {\cal O}x_1 \oplus {\cal O}x_2\oplus
{\cal O}x_3$ where $\{x_1, x_2, x_3\}$ is some basis for $V$. Thus the columns of a 
non-singular $3\times3$ matrix with entries in ${\bf Z}_p(t)$ defines a lattice. 
The standard lattice is the one
corresponding to the identity matrix. We define two lattices $L$  and $L'$ to be 
equivalent, if for some $\lambda \in {\bf Z}_p(t)^*$ we have $L = \lambda L'$.
 We denote equivalence class by $[L]$. The building $\Delta$ is defined as   
a flag complex in the following way. The points are
equivalence classes of lattices, and $[L_0], .... , [L_k]$ span  a $k$--simplex 
(in our situation $k = 0 , 1 , 2$ are the  only possibilities) if and only if one 
can find representatives so that 
$\pi L_0 \subset L_1 \subset ... \subset L_k \subset L_0$. 

All 2--simplices   are of the form $\{ [ x_1, x_2 , x_3 ]
, [ x_1, x_2, \pi x_3 ] , [ x_1, \pi x_2, \pi x_3 ] \} $; this is usually 
referred to as a {\em chamber} and denoted by $C$. Clearly the group 
 $GL_3({\bf Z}_p(t))$ acts on lattices and one sees easily that incidence 
 is preserved, so that the group acts simplicially on
$\Delta$.  It  is shown in \cite{Br} that this  building is a so-called {\em
Euclidean building}, in particular, it is contractible and can be equipped with a metric
which makes it into a $CAT(0)$  space and for which $GL_3({\bf Z}_p(t))$ acts 
as a group of isometries. The metric is
such that each $2$ dimensional simplex is isometric to a unit Euclidean
equilateral triangle. 

We now return to our situation. One of the difficulties of dealing with 
representations of braid groups is that it
is extremely difficult to determine exactly which matrices are in the image. We bypass 
this by dealing with a group which contains $im(\beta_4\otimes{\bf Z}_p)$. To define
this group, we recall that it was shown by Squier \cite{Sq} that the  Burau representation is 
unitary in the sense that there is a matrix  
\[ J = \left( \begin{array}{ccc} -(s+1/s) & 1/s & 0 \\
 s &  -(s+1/s) & 1/s \\
 0 &  s & -(s+1/s)  \end{array} \right)  \]
with the property that $A^*JA = J$ for all $A \in im(\beta_n)$. Here the involution $*$ comes from
extending the involution of ${\bf Z}_p[t,t^{-1}]$ generated by $t \rightarrow 1/t$
to the matrix group by  $(a_{i,j})^* = ( a_{j,i}^*)$, where $s^2 = t$. 

We define  the subgroup  $Isom_J(\Delta)$ of $GL( 3 , {\bf Z}_p(t))$ to be those 
matrices with Laurent polynomial entries which are unitary for the form $J$. The advantage of
dealing with this subgroup is that the condition that a matrix lies inside $Isom_J(\Delta)$
is easily used. 

The strategy now is to examine the action of $Isom_J(\Delta)$ on $\Delta$. This 
is interesting in its own right. Moreover, the greater ease of dealing with this
subgroup means that we are able to  compute the complex $\Delta/Isom_J(\Delta)$ 
together with all vertex, edge and $2$--simplex 
stabilisers. We then appeal to results of Haefliger \cite{H} to compute a
presentation for the group $Isom_J(\Delta)$. 

Now recall that homotheties act trivially on $\Delta$ so that the presentation for
$Isom_J(\Delta)$ is to be compared with the following presentation of $B_4/centre(B_4)$:
\begin{lemma}
The group $B_4/centre(B_4)$ is presented as
\[ \langle\, x , y \;\;|\;\;x^4 = y^3 = 1 \quad [ x^2 , yxy ] = 1\, \rangle \]
where $ x = \sigma_1\sigma_2\sigma_3$ and $y = x\sigma_1$.
\end{lemma}
This is presumably well known to the experts---it is derived in \cite{CL}. The starting point 
for this work is:
\begin{lemma}
\label{stab_group_points}
The group $stab_J(I)$ acts on $\Delta$ as a finite group.
\end{lemma}
\begin{proof}[Sketch of proof] 
If $A \in stab_J(I)$, then its action on $\Delta$ is
unchanged by homothety and it's easily seen that we can adjust any such
$A$ by applying $\pm t^k$ so that $A \in SL(3 , {\cal O})$. Rewriting the unitary
condition as $A^* = JA^{-1}J^{-1}$ and noting that $J \in GL(3 , {\cal O})$, we 
see that $A^* \in SL(3 , {\cal O})$. However the only matrices with Laurent polynomial
entries for which $A$ and $A^*$ have all entries valuing positively are the 
constant matrices. 

Thus we have shown that the only such $A$ have constant entries up to homothety. In
particular, they are unchanged by setting $t = 1$, so that  $stab_J(I)$ can be 
regarded as a subgroup of the finite group $GL(3,{\bf Z}_p)$, completing the 
proof.  
\end{proof}

This has the immediate corollary:
\begin{corollary}
For every vertex $v \in \Delta$, $stab_J(v)$ is a finite group.
\end{corollary}
\begin{proof} The building $\Delta$ is locally finite, in fact the link of every
vertex is the flag manifold in the vector space ${\bf Z}_p^3$. The stabiliser of
any vertex acts on this set as a group of permutations, so by passing to a subgroup
of finite index in $stab_J(v)$ we obtain a subgroup which acts as the identity on
all vertices in the link. Since every vertex is connected to $I$ by some chain of
vertices, we see that for every $v$, there is a subgroup of finite index which lies
inside $stab_J(I)$, a finite group.  \end{proof}

We now focus on the case $p = 3$.  In this case one finds by calculation:
\begin{theorem}
\label{I_stab}
At the prime $3$,  group $stab_J(I)$ acts on $\Delta$ as ${\bf Z}_4 \cong \langle x \rangle$.
\end{theorem}
\begin{remark}
For $p = 2, 3 ,5$, the group $stab_J(I)$ acts as the cyclic group ${\bf Z}_4$. For $p = 7$ it is cyclic of order
$8$ and for $p = 11$, cyclic of order $12$.
\end{remark}
One important difference between the case  $ p = 2$ and that of the larger primes is that it is one of 
the consequences of  the results  of \cite{CL} that  $Isom_J(\Delta(2)) \cong  im(\beta_4\otimes{\bf Z}_2)$, 
this is not so for (at least some and conjecturally all) primes $p \geq 3$. In particular, for $p = 3$, we are able to
construct (see below) an element $u \in Isom_J(\Delta(3))$ which has order $6$; 
it is easy to see that this element does not lie in the subgroup $im(\beta_4\otimes{\bf Z}_3)$.   
Its matrix is given by:
\[ u = \left( \begin{array}{ccc} 2 + t + t^2 & 2 + t^2 & 2 + 2 t + 2 t^2 \\
 2 + 2 t^2 & 2 + t + 2 t^2 & 2 + t + t^2 \\
  2 + t &  2 + t & 2 + 2 t  \end{array} \right)  \]
However, having noted this difference, the qualitative picture of the quotient complex 
is very similar to the case $ p = 2$; the complex consists of a compact piece coming
from behaviour of groups close to the identity lattice, together with a single annular
end. Application of Haefliger's methods yields the following group theoretic result:
 \begin{theorem}
 When $p = 3$, the group $Isom_J(\Delta)$ is presented as:\\
 Generators: $\quad  x , y , u $\\
 Relations:
 \begin{enumerate} \def\labelenumi{{\rm(\theenumi)}}
 \item $x^4 = y^3 = u^6 = 1 $
 \item $[x^2 , yxy] = 1$
 \item $[x , u^{-1}x^{-1}y^{-1}xyxy ] = 1 $
 \item $[yxy , u^{-1}x^{-1}y^{-1}xyxy ] = 1 $
 \item $[xyx , u^2] = 1$
 \item $[x^2yx , u^3] = 1$
 \item $(u^2x^2yx)^2 = (x^2yxu^2)^2$
 \item \label{other_relations} Infinitely many other relations to do with nilpotence.
 
 \end{enumerate}
 \end{theorem}
 Of course the verification that these relations hold is a trivial matter of multiplying
 matrices modulo $3$. We remark that the relations contained in (\ref{other_relations}) are explicitally
known.

We claim that a computer application of the Reidemeister--Schreier algorithm contained in the computer program
{\em GAP} applied to the presentation involving the first seven relations proves:
 \begin{corollary}
 The index $[Isom_J(\Delta) : \langle x , y \rangle ]$ is finite.
 \end{corollary}
This index is a divisor of $162$. 
The corollary already implies that $im(\beta_4\otimes{\bf Z}_3)$ is not faithful. One way to see this
 is that one sees easily (for example from the matrix representation) that the element 
 $w = u^{-1}x^{-1}y^{-1}xyxy$ has infinite order. The presentation implies that it commutes with 
 $x$. However, since $[Isom_J(\Delta) : \langle x , y \rangle ]$ is finite, some power of $w$ lies in the
 subgroup generated by $x$ and $y$ and this gives an unexpected element commuting with $x$.
 Alternatively, in the course of the proof, one discovers that $Isom_J(\Delta)$ contains
 arbitrarily large soluble subgroups and this can also be used to show that the representation
 is not faithful. In fact, one can be more specific; the computer can be used to give a 
 presentation for the subgroup generated  by $\langle x , y \rangle$;  one finds for example, 
 that there is a relation (where $\xb=x^{-1}$ and $\yb=y^{-1}$):
\[
\begin{array}{c}
  \xb. \yb. \xb. y. x. \yb. x. y. x. \yb. x. \yb. \xb. y. x. y. \xb. \yb. \xb. \yb. 
      \xb. y. x. \yb. \xb. \yb. \xb. y. x. \yb. x. y. x. \yb. x. \yb.   \\

\xb. y. x. y. x. 
      \yb. \xb. \yb. x. y. x. \yb. \xb. \yb. \xb. y. x. \yb. x. y. x. \yb. x. \yb. \xb. 
      y. x. y. x. \yb. \xb. \yb. x. y. x. \yb   
\end{array}
\]         
That this relation does not hold in the braid group is easily checked by computing the integral Burau matrix.

 \section{Outline of the proof for $p = 3$}
 In spirit, if not in detail, the proof follows the ideas introduced in \cite{CL},  to which we 
refer the reader. We work outwards from the identity lattice, successively identifying point
stabilisers. This enables us to find representatives for each orbit and hence build the quotient
complex. The compact part alluded to above comes from the action of the group on vertices fairly
close to the orbit of the identity; as one moves farther away there is a certain amount  of stabilisation
and it is this which gives rise to the single annular end.

We refer to the orbit of the lattice $I$ as the {\em group points}. The result  Lemma
\ref{stab_group_points}   shows that every 
group point has stabiliser ${\bf Z}_4$. We recall the link of any vertex may be considered as the flag geometry
of the vector space ${\bf Z}_3^3$, so that every vertex has $26$ points in its link, and each vertex in the link
is adjacent to four other vertices in the link.

We need to recall the notation introduced in \cite{CL}. We make a (noncanonical) choice of 
representative lattices for each  of the $26$ vertices by writing down matrices whose
columns define the lattice. Subsequent vertices are coded by using these matrices, regarded
as elements of $GL(3,{\bf Z}_3[t,t^{-1}])$ as acting on $\Delta$. As an example, denoting
the matrix representative chosen for the thirteenth vertex by $M_{13}$, then the representative
elements in the link of the the thirteenth vertex are chosen to be 
$M_{13}.M_j$ for $1 \leq j \leq 26$. Of course, one vertex has several names in this notation,
for example the identity vertex appears in the link of each of its vertices.

The first task is to examine how many group points lie in the link of the identity. 
 \begin{lemma}
 $Link(I)$ contains precisely $18$ group points: 
\[
\begin{array}{c}
y , \quad y^2 , \quad x.y ,\quad x.y^2 , \quad  x^2.y ,\quad  x^2.y^2,\\
\quad x^3.y , \quad x^3.y^2 , \quad y.x.y , \quad (y.x.y)^{-1}, \quad x.y.x.y, \quad x.(y.x.y)^{-1},\\ 
\quad w ,\quad  w^{-1} ,\quad (y.x.y)^{-1}.w ,\quad  x.(y.x.y)^{-1}.w ,\quad  y.x.y.w^{-1} ,
\quad x.y.x.y.w^{-1},
\end{array}
\]
 where $w$ is the element introduced at the end of section 1.
\end{lemma}
Of course the fact that these are all group points is immediate and the fact that
they are distance one from $I$ is a calculation. The content of the lemma is
that there are no more group points. This proved by noting that the lattice 
\[ M_{19} = \left( \begin{array}{ccc} 1 & 0 & t\\
  0 & t  &  0 \\
  0  &  0 &  t  \end{array} \right)  \]
is in the link of the identity and is stabilised by the element $u$. Thus it cannot
be a group point as its stabiliser contains an element of order $6$. The action of
known group elements now accounts for all the other elements in $Link(I)$. 

We indicate briefly how one can construct any isometries which may exist in the stabiliser of $M_{19}$,
in particular, how one can find the element $u$. This involves an elaboration of the method used in
Lemma \ref{stab_group_points}. 

Suppose that $g \in Isom_J(\Delta)$ has $g[M_{19}] = [M_{19}]$. The definition shows that
this is the same as the existence of an element $\alpha \in GL_3({\cal O})$ with $g.M_{19} = M_{19}.\alpha$.
The form of the elements $M_{19}$ and $g$ means that $\alpha$ has Laurent polynomial entries. Then
\[ \alpha^*(M_{19}^*JM_{19})\alpha = (M_{19}\alpha)^*.J.(M_{19}\alpha) = (g.M_{19})^*.J.(g.M_{19}) = M_{19}^*JM_{19} \]
since $g$ is an isometry. It follows that $\alpha$ is an isometry of the form $M_{19}^*JM_{19}$ and although unlike
Lemma \ref{stab_group_points}, this form does not have its matrix lying in $GL_3({\cal O})$, we have a bound on the valuations of
its entries, so that exactly as in the lemma, we have a bound on the valuations possible for the entries of $\alpha$. Since
we are dealing with a fixed finite field, it follows that there are only a finite number of possibilities for the entries of 
$\alpha$ and one can check by direct enumeration which of these make $M_{19}\alpha M_{19}^{-1}$ into a $J$ isometry. 
(In fact sharper, more practical methods exist, but this would  take us too far afield.) 

We now give some indication of how one can give complete descriptions of all vertex stabilisers. The idea is to
work outwards from the identity; it turns out that we need no more elements than those we have already introduced.

Recalling the notation defined above, a calculation shows that that action of $u$ on its link is given by the permutation
\[ 
\begin{array}{c}
(7^*)(11^*)(18^*)(23^*)(3^* 13^*)(6^* 8^*)(14^* 24^* 26^*)(17^* 21^* 20^*)\\
(1^* 5^* 12^* 22^* 19^* 10^*)(4^* 16^* 25^* 15^* 2^* 9^*) 
\end{array}
\]
where $x^*$ is shorthand for $M_{19}.x$. The two six cycles consist of $12$  group points, 
($I = 2^*$), there are $14$ points in the orbit of $M_{19}$ and two remaining, as yet unidentified points, 
$7^*$ and $11^*$. Points in the orbit of $M_{19}$ we refer to as $n$--points. Observe that neither of the unidentified 
points can be group points as they contain an element of order $6$ in their stabiliser.

Using this information we now show:
\begin{lemma}
\label{m19_stab}
The group $stab_J(M_{19})$ acts on $\Delta$ as a finite group ${\bf Z}_6 \cong \langle u \rangle$.
\end{lemma}

\begin{proof}[Sketch of Proof] First consider the map $i_0 \co stab_J(M_{19}) \rightarrow Aut(Link(M_{19}))$. We begin
by noting that this map is injective, for any element of the kernel must fix every vertex in $Link(M_{19})$,
in particular the vertex $I$, so that the kernel can only consist of powers of the element $x$. However, one
checks that no element of the group $\langle x \rangle$ other than the identity fixes $M_{19}$ proving the assertion.

We refer to the above permutation, where we recall the vertex $2^*$ is the
identity vertex. Pick an element $\gamma \in stab_J(M_{19})$; it is type-preserving so that it must map the group
points in $Link(M_{19})$ which correspond to lines back to lines, and those which correspond to planes to planes.
Since $u$ acts transitively on this orbit, we can find some power of $u$ so that $u^k.\gamma$ fixes the vertex
$2^*$. Now exactly as in the previous paragraph, we deduce that $u^k.\gamma = I$, so that $\gamma$ is a power
of $u$ as required.   
\end{proof}

We now analyse the two new points $7^*$ and $11^*$. We have already shown that these are not group points; we now
show that they are not $n$--points.

Firstly, one finds that $xyx(7^*) = 11^*$, so that this is only one orbit of point and moreover  that $u$ acts 
as an element of order two on $Link(7^*)$. Moreover, we can construct a 
potentially new element in $stab_J(7^*)$ namely $u_1 = (xyx)^{-1}.u.xyx$. A calculation reveals that the 
action of the group $\langle u , u_1\rangle$ on $Link(7^*)$ is the dihedral group $D_3$. It now follows from \ref{I_stab} and
\ref{m19_stab} that the orbit of $7^*$ is distinct from that of the group and $n$--points.

In fact, the stabiliser is larger than this and one finds that there
is an element $h \in \langle x,y,u\rangle $ of order $3$ 
which commutes with this dihedral group. 
\[  h  = \left( \begin{array}{ccc} 1 + t^4  & 1 + t^2 + t^4 & 1 + t + 2t^2 + 2t^3\\
  2t + 2t^2 + 2t^4 & 2 + t^2 + 2t^4 &  2 + 2t + t^2 + t^3 \\
  0  &  0 &  2t^2  \end{array} \right)  \]
We omit the arguments which identify the stabilisers of these two points, 
as this is slightly special, however the results are that one shows successively:
\begin{lemma}
\label{seven_star}
The map $i_1 \co stab_J(7^*) \rightarrow Aut(Link(7^*))$ has $ im(i_1) \cong {\bf Z}_3 \times D_3 $. 
\end{lemma}
\begin{corollary}
The group $stab_J(7^*)$ has order $54$ with structure given by the nonsplit central extension:
$$1 \rightarrow \langle u^2 \rangle \cong {\bf Z}_3 \rightarrow stab_J(7^*) \rightarrow {\bf Z}_3 \times D_3 \rightarrow 1$$
\end{corollary}
The orbit type for the action of $stab_J(7^*)$ acting on its stabiliser is $\{9,9,3,3,1,1\}$ where the orbits of
size $9$ are $n$--points, the orbits of size $3$ are of type $7^*$ and there are two points as yet unaccounted
for, namely $M_{19}.M_7.M_7$ and $M_{19}.M_7.M_{11}$ for  which we adopt the notational shorthand
$7^{(2)}$ and $11^{(2)}$. As above, $xyx(7^{(2)}) = 11^{(2)}$.

This is the point at which the behaviour stabilises. For later use, it is more convenient to define for $i \geq 0$, 
a sequence of elements $\alpha_{i+1} = (xyx)^{-i}u.u_1(xyx)^i$. Then we have:
\begin{theorem}
\label{seven_stabs}
For $k \geq 2$, the map $i_k \co stab_J(7^{(k)}) \rightarrow Aut(Link(7^{(k)}))$ has image of order $54$.

Moreover, $stab_J(7^{(k)})$ is generated by the elements $u$, $h$, $\alpha_1$, .... , $\alpha_k$.
\end{theorem}

\begin{proof}[Sketch Proof] The argument is inductive; we explain the step $k = 2$ which contains all the essential ingredients.
We set $H(2) = \langle u , \alpha_1 , \alpha_2 \rangle \leq stab_J(7^{(2)})$. Note that every element of $H(2)$ stabilises 
$7^{(3)}$ and $11^{(3)}$. We refer to Figure 1, which shows the hexagon $Link(7^{(2)})/H(2)$. 
Our  claim is that no element of  $\eta \in stab_J(7^{(2)})$  can move $7^{(3)}$.

\begin{figure}[htb]\cl{%
\epsfysize=180pt                        
\epsfbox{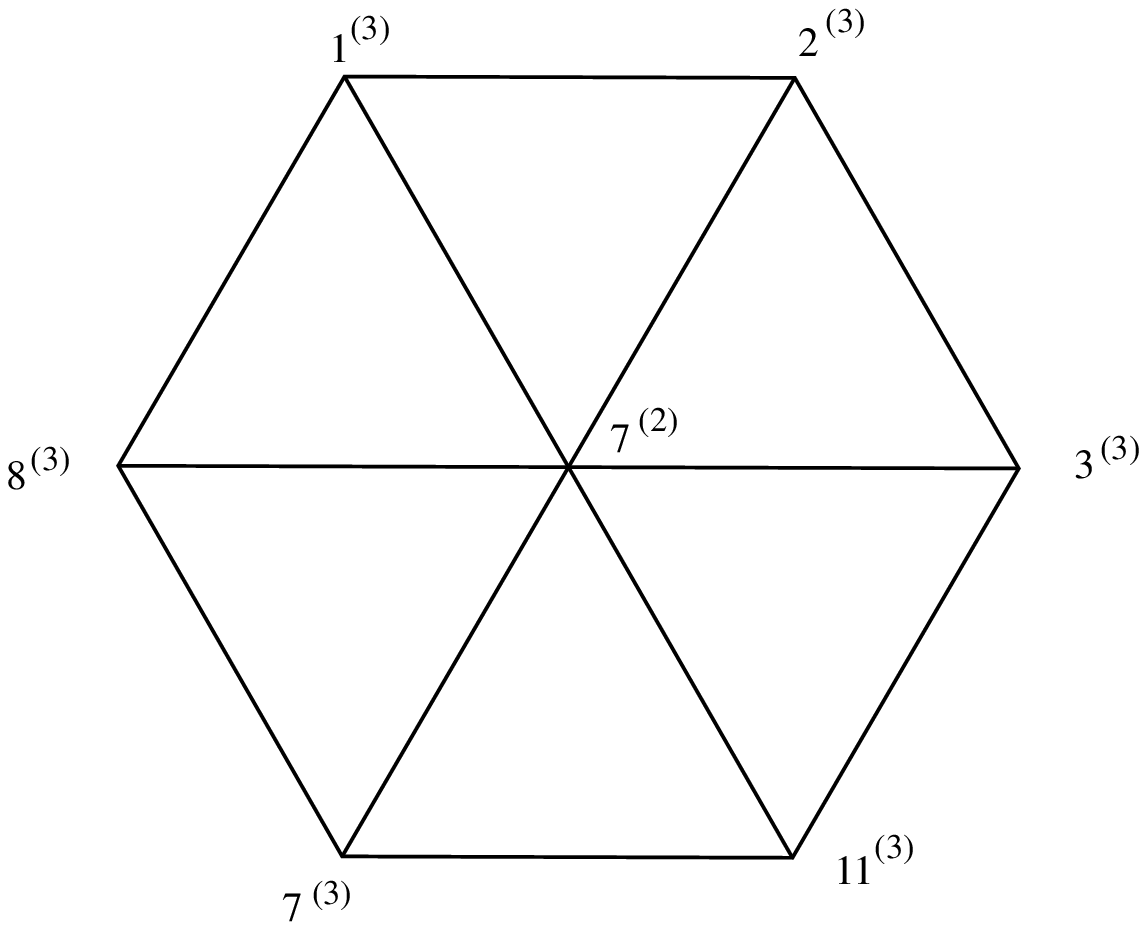}} 
\cl{\small Figure 1}
\end{figure}

We argue as follows. Note that since elements in vertex stablisers are type preserving, the only
possibilities for $\eta(7^{(3)})$ (assuming that it is moved) are the $H(2)$ orbit of $1^{(3)}$ or the 
$H(2)$ orbit of $3^{(3)}$.

However, the former orbit contains $9$ elements and the latter $3$, so that in any case, $\eta$ must move some element
in the $H(2)$ orbit of $1^{(3)}$ back into this orbit. By composing with an element of $H(2)$, we see that this implies
the existence of an element moving $7^{(3)}$ lying in $stab_J(k^{(3)}) \cap im(i_2)$ where $k^{(3)}$ lies in the $H(2)$ 
orbit of $1^{(3)}$. After conjugating by an element of $H(2)$, we may assume that this element lies in $stab_J(1^{(3)}) 
= (xyx)^{-1}stab_J(2^{(3)})xyx$. But $2^{(3)} = 7^*$, so that $stab_J(1^{(3)}) = (xyx)^{-1}stab_J(7^*)xyx$. 
An examination of the generating elements shows that no element of this latter group moves $7^{(3)}$, a contradiction.

A similar argument establishes that $stab_J(7^{(2)})$ stabilises $11^{(3)}$.

We now show that $im(i_2)$ is a group of order at most $54$. The reason is this: All of $im(i_2)$ stabilises $7^{(3)}$
hence permutes the four points in the link adjacent to it, however one of these points is $11^{(3)}$, which is 
also fixed by the whole group. Therefore by passing to a subgroup of $im(i_2)$ of index at most $3$ we stabilise
the point $3^{(3)}$. Arguing similarly for $3^{(3)}$, we deduce that $im(i_2)$ contains a subgroup of index at
most $9$ which stabilises $2^{(3)} = 7^*$. This is a group whose structure is already completely determined and
one finds that $stab_J(7^*)$ acts on $Link(7^{(2)})$ as a group of order $6$, proving the claim.

Now the group $H(2)$ is easily analysed; in particular, one shows easily that it acts on the link as a
group of order $54$. This establishes that $i_2(H(2)) = im(i_2)$ as required. 

The kernel of the map $i_2 \co stab_J(7^{(2)}) \rightarrow Aut(Link(7^{(2)}))$  is a subgroup of $stab_J(7^*)$. 
Recalling that $H(2)$ is generated by $u$, $\alpha_1$ and $\alpha_2$, it follows that $u$, $h$, $\alpha_1$ and 
$\alpha_2$ generate $stab(7^{**})$, completing the first step of the  induction. 
\end{proof}

Given this theorem, one can now give a complete description of the groups $stab_J(7^{(k)})$ by analysing
how $ker(i_k) \leq stab_J(7^{(k-1)})$ acts on $Link(7^{(k)})$. As a consequence, one proves that the 
group $stab_J(7^{(k)})$ has order $2.3^{2k+1}$. It follows immediately that the orbits $7^{(k)}$ are all 
distinct.

We recap our progress so far. From the information that the stabiliser of the lattice $I$ is a cyclic group
of order four, we have identified the stabiliser of every vertex in the building; this information suffices to
deduce that the orbits for the action of $Isom_J(\Delta)$ on $\Delta$ are precisely $I , M_{19} , 7^{(k)}$ for
$k \geq 1$. Moreover, this already shows:
\begin{corollary}
The group $Isom_J(\Delta)$ is generated by $x$, $y$ and $u$.
\end{corollary}
The construction of the entire complex $\Delta/Isom_J(\Delta)$ rests largely on the work set forth above and
we shall not go into it in detail. Broadly it involves two steps: The identification of a candidate set of
orbits of edges and triangles coming from the action of fairly short elements, followed by the proof that
no further identifications are possible. This latter step is accomplished by the detailed understanding 
we have developed of the vertex stabilisers. This task gets easier as one moves further away from the
group points, as stabilisers get larger and there are less orbits to be considered; eventually the action
of stabilisers on links becomes constant. As a result, the complex has a fairly natural decomposition into
two pieces; a compact part and some ``tubes". We refer the reader to \cite{CL} for details in the case
$p = 2$. For example, a picture of the tube comes from the concatenation of hexagons shown in Figure 2.

\begin{figure}[htb]\cl{%
\epsfysize=200pt                        
 \epsfbox{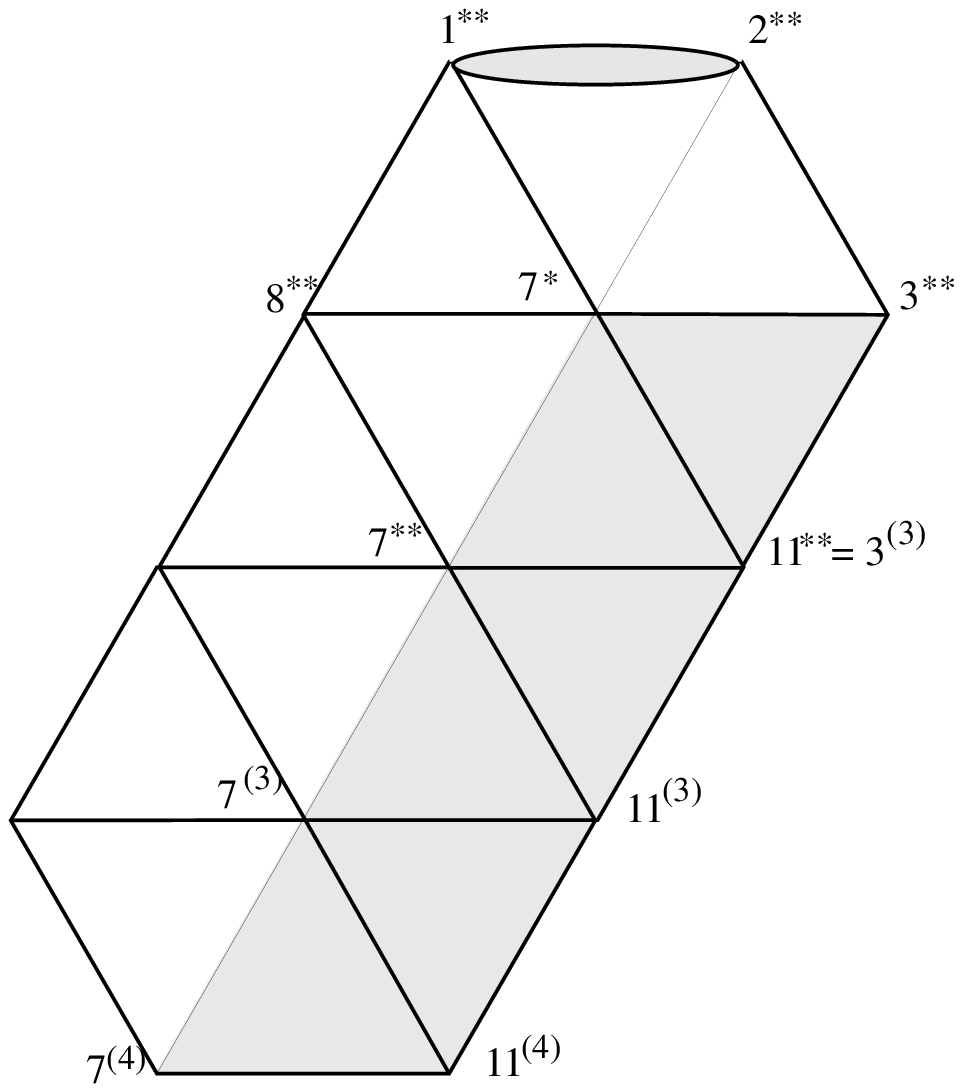}} 
\vspace{2mm}
\cl{\small Figure 2}
\end{figure}

\section{The case $p = 5$}
The analysis in this case follows the same outline as indicated above, though of course the details become
much more complicated. Nonetheless, one obtains a presentation of the group $Isom_J(\Delta(5))$. The
quotient complex has interesting features not present in the first two cases; for example in contrast to the 
cases $p = 2$ and $p = 3$, the complex  which emerges has three annular ends.

Once again one finds extra elements in $Isom_J(\Delta(5))$ which it turns out do not lie in the group generated by
$x$ and $y$. The simplest of these is the element $\beta_2$ shown below:
\[ \beta_2 = \left( \begin{array}{ccc} 4 & 1 + 2t + 2t^2 & 3 + t\\
  1 + t & 4 + 2t  &  2 + 2t \\
  1  &  4 + 3t + 4t^2 &  2 + 2t  \end{array} \right)  \]
This is an element of order $4$ and one finds that:
\begin{theorem}
The group $Isom_J(\Delta(5))$ is generated by $x$, $y$ and $\beta_2$
\end{theorem}
In fact, we are able to complete all the analysis up until the very last step and in particular, we are able to find
a presentation of the group $Isom_J(\Delta(5))$. It is rather complicated and {\em GAP} was unable to show that
the index $[Isom_J(\Delta(5)) : \langle x , y \rangle ]$ was finite. We have been unable to prove that it is infinite and unable
to analyse the situation sufficiently to prove or disprove that $\langle x , y \rangle$ contains no extra relations.

\ppar

{\bf Acknowledgment}\qua Both authors are supported in part by the NSF.

 \vspace{0.2in}{\small\sl\parskip0pt
Department of Mathematics\\
University of California\\
Santa Barbara, CA 93106, USA\par
\medskip
\rm Email:\stdspace\tt cooper@math.ucsb.edu, long@math.ucsb.edu\par}
\recd 
 
\end{document}